\documentclass[11pt,twoside]{article}
\usepackage{amssymb}

{\catcode`\@=11 \gdef\ps@myheadings{\let\@mkboth\@gobbletwo
 \def\@oddhead
          {\vbox{\noindent
                {\small SUPREMUM OF A RANDOM WALK}
                                         \hfill\rm\thepage\vskip 4pt}}%
 \def\@oddfoot{}
 \def\@evenhead
          {\vbox{\noindent\rm\thepage\hfill\small
          D. DENISOV, S. FOSS AND D. KORSHUNOV\vskip 4pt}}%
 \def\@evenfoot{}\def\chaptermark##1{}\def\sectionmark##1{}%
 \def\subsectionmark##1{}}

\gdef\@begintheorem#1#2{\trivlist \item[\hskip
\labelsep{\bf #1\ #2.}]}
\gdef\@opargbegintheorem#1#2#3{\trivlist
      \item[\hskip \labelsep{\bf #1\ #2\ {\rm (#3)}.}]}
\gdef\@endtheorem{\endtrivlist}
 }

\font\scal=rsfs10 at 11pt
\font\scall=rsfs10
\newcounter{remark}
\newcommand\remark
       {\refstepcounter{remark}
                  {R\,e\,m\,a\,r\,k~\arabic{remark}.\ }}
\newcommand\proof{\noindent{\it Proof}}
\newcounter{definition}
\newcommand\definition[1]
       {\refstepcounter{definition}
                  {\noindent\bf Definition~\arabic{definition}{\rm#1}.\ }}

\newlength{\myskip}
\newtheorem{Theorem}{Theorem}
\newtheorem{Proposition}{Proposition}
\newtheorem{Lemma}{Lemma}
\newtheorem{Corollary}{Corollary}
\newtheorem{Example}{Example}

\newcommand\mysection[1]{
              \refstepcounter{section}
              \section*{\normalsize\bf\thesection.~#1}
                         }

\begin{document}
\thispagestyle{empty}

\section*{\Large\bf
           Tail Asymptotics for the Supremum of a Random Walk
           when the Mean Is Not Finite\footnote{Supported
       by EPSRC grant No.~R58765/01, INTAS grant No.~265
              and RFBR grant No.~02-01-00358}}
          
\section*{\normalsize\rm
DENIS DENISOV \hfill denisov@ma.hw.ac.uk\\
SERGUEI FOSS \hfill foss@ma.hw.ac.uk\\
DIMA KORSHUNOV \hfill korshunov@ma.hw.ac.uk\\
{\em Department of Actuarial Mathematics and Statistics,
School of Mathematical and Computer Siences,
Heriot-Watt University,
Edinburgh EH14 4AS, Scotland}}

\vspace{5mm}

{\small {\bf Abstract.}
We consider the sums $S_n=\xi_1+\cdots+\xi_n$ of
independent identically distributed random variables.
We do not assume that the $\xi$'s have a finite mean.
Under subexponential type conditions on distribution
of the summands, we find the asymptotics of the probability
${\bf P}\{M>x\}$ as $x\to\infty$,
provided that   $M=\sup\{S_n,\ n\ge1\}$ is a proper random variable.
Special attention is paid to the case of tails which are
regularly varying at infinity.

We provide some sufficient conditions for the integrated
weighted tail distribution to be subexponential. We supplement these
conditions by a number of examples which cover both
the infinite- and the finite-mean cases.
In particular, we show that the subexponentiality of 
distribution $F$ does not imply the subexponentiality 
of its integrated tail distribution $F^I$.

\vspace{5mm}

\noindent{\bf Keywords}: supremum of sums of random variables,
large deviation probabilities,
subexponential distribution,
integrated weighted tail distribution
}

\mysection{Introduction}

Let $\xi$, $\xi_1$, $\xi_2$, \ldots\ be independent
random variables with common non-degenerate distribution~$F$
on the real line ${\bf R}$.
We let $F(x)=F((-\infty,x])$ and $\overline F(x)=1-F(x)$.
In general, for any distribution $G$, we denote its tail by
$\overline G(x)=G((x,\infty))$.
In this paper, an important role is played by
{\it the negative truncated mean function}
\begin{eqnarray*}
m(x) &\equiv& {\bf E}\min\{\xi^-,x\}
=\int_0^x {\bf P}\{\xi^->y\}\,dy,\quad x\ge0,
\end{eqnarray*}
where $\xi^-=\max\{-\xi,0\}$; the function $m(x)$ is continuous,
$m(0)=0$ and $m(x)>0$ for any $x>0$.

Put $S_0=0$, $S_n=\xi_1+\cdots+\xi_n$, and
$$
M=\sup\ \{S_n,\ n\ge 0\}.
$$
Our main assumption is that $M$ is finite a.s.
The latter occurs if and only if $S_n\to-\infty$ as $n\to\infty$
with probability one (see Theorem 1 in
[\ref{F}, Chapter XII, Section 2]). It is known that

(i) if ${\bf E}|\xi|<\infty$, then $S_n\to-\infty$ a.s.
as $n\to\infty$ if and only if ${\bf E}\xi<0$;

(ii) if ${\bf E}|\xi|=\infty$, then $S_n\to-\infty$
a.s. as $n\to\infty$ if and only if
\begin{equation}\label{finite.M}
\int_0^\infty\frac{x}{m(x)}\,F(dx)\ \mbox{ is finite,}
\end{equation}
see Corollary 1 in [\ref{Erickson}]. Note that
the function $\frac{x}{m(x)}$ is increasing, since
\begin{eqnarray}\label{derivative.of.t.mt}
\frac{d}{dx}\frac{x}{m(x)}
&=& \frac{m(x)-xm'(x)}{m^2(x)}
=\frac{m(x)-x{\bf P}\{\xi^->x\}}{m^2(x)}\ge 0.
\end{eqnarray}

In the case (ii), $m(x)\to\infty$ as $x\to\infty$, with necessity.
Roughly speaking, the condition (\ref{finite.M}) means that
the right tail of the distribution $F$ is lighter than the left one.

The main goal of the present paper is to investigate
the asymptotic behaviour of the probability ${\bf P}\{M>x\}$
as $x\to\infty$ when the distribution of the summands is heavy-tailed.
As far as applications are concerned, (a) in queueing,
$M$ coincides in distribution with the stationary waiting time in the corresponding
$GI/G/1$ queue; (b) in risk theory, ${\bf P}\{M>x\}$
is the probability of ruin.

We recall the definitions of some classes of functions
and distributions which will be used in the sequel.

\definition{}
The function $f$ is called {\em long-tailed}
if, for any fixed~$t$, the limit of the ratio
$f(x+t)/f(x)$ is equal to~$1$ as $x\to\infty$.
We say that the distribution~$G$ is long-tailed
(and write $G\in\mbox{\scal L}\,$)
if the function~$\overline G(x)$ is long-tailed.

\definition{}
The distribution~$G$ on ${\bf R}^+$
with unbounded support belongs to the class
$\mbox{\scal S}\,$ (and is called {\it a subexponential distribution})
if the convolution tail~$\overline{G{*}G}(x)$
is asymptotically equivalent to~$2\overline G(x)$ as $x\to\infty$.

It is shown in [\ref{Ch}] that any subexponential
distribution~$G$ is long-tailed with necessity.
Sufficient conditions for some distribution to
belong to the class~$\mbox{\scal S}\,$ may be found,
for example, in [\ref{Ch}, \ref{Kl}, \ref{T}].
The class~$\mbox{\scal S}\,$ includes, in particular,
the following distributions on $[0,\infty)$:
(i)~any distribution $G$ whose tail $\overline G(x)$ is regularly varying
at infinity with index $\alpha<0$, that is, for any fixed $t>0$,
$\overline G(xt)\sim t^\alpha\overline G(x)$ as $x\to\infty$;
(ii)~the lognormal distribution with the density
$e^{-(\ln x-\ln\alpha)^2/2\sigma^2}/x\sqrt{2\pi\sigma^2}$
with ${\alpha>0}$;
(iii)~the Weibull distribution with the tail
$\overline G(x)=e^{-x^\alpha}$ with ${\alpha\in(0,1)}$.

It is known (see [\ref{V}] and [\ref{EV}]) that, if ${\bf E}\xi=-a$ is finite
negative number and the integrated tail distribution~$F^I$,
\begin{equation}\label{def_G_I}
\overline{F^I}(x) =
\min\Big(1,\ \int_0^\infty \overline F(x+u)\,du \Big),
\quad x>0,
\end{equation}
is subexponential, then the distribution tail of the maximum
of sums is equivalent, up to a constant,
to the integrated tail of the distribution of one summand, that is,
\begin{equation}\label{as_infty}
{\bf P}\{M>x\} \sim \overline{F^I}(x)/a\quad \mbox{as } x\to\infty.
\end{equation}
The converse is also true (see [\ref{K}]):
if the asymptotic (\ref{as_infty})
holds, then the integrated tail distribution $F^I$
is subexponential.

In the present paper, we consider mainly the case where the $\xi$'s
have infinite mean. In this case, we should assume
${\bf E}\xi^-=\infty$, otherwise $M=\infty$.
Without further assumptions, we can provide lower and upper bounds only.

\begin{Theorem}\label{bounds.gen}
Suppose ${\bf E}\xi^-=\infty$ and the condition {\rm(\ref{finite.M})} holds.
Let the distribution $F$ be long-tailed and the distribution $G_1$ with the tail
\begin{eqnarray}\label{tail.G.1}
\overline G_1(x) &=&
\min\Biggl(1,\int_0^\infty\overline F(x+t)\,d\frac{t}{m(t)}\Biggr)
\end{eqnarray}
be subexponential. Then the following estimates hold:
\begin{eqnarray*}
1 \le \liminf_{x\to\infty}\frac{{\bf P}\{M>x\}}{\overline G_1(x)}
&\le& \limsup_{x\to\infty}\frac{{\bf P}\{M>x\}}{\overline G_1(x)} \le 2.
\end{eqnarray*}
\end{Theorem}

In the case where the function $m(x)$ is regularly varying,
we get the following sharp asymptotics
(the symbol $\Gamma$ stands for the {\it Gamma function}):

\begin{Theorem}\label{asy.reg.-}
Suppose ${\bf E}\xi^-=\infty$ and the condition {\rm(\ref{finite.M})} holds.
Let $m(x)$ be regularly varying at infinity with index
$1-\alpha\in[0,1]$. If the distribution $F$ is long-tailed and
the distribution $G_1$ with the tail {\rm(\ref{tail.G.1})}
is subexponential, then
\begin{eqnarray}\label{M.asy.reg.1}
{\bf P}\{M>x\} &\sim&
\frac{\overline G_1(x)}{\Gamma(1+\alpha)\Gamma(2-\alpha)}
\quad\mbox{ as } x\to\infty.
\end{eqnarray}
If $\alpha\in(0,1]$, then the assumption of the subexponentiality
of $G_1$ can be replaced by that of the subexponentiality of the
distribution $G_2$ with tail
\begin{eqnarray}\label{tail.G.2}
\overline G_2(x) &=&
\min\Biggl(1,\int_1^\infty\frac{\overline F(x+t)}{m(t)}\,dt\Biggr),
\end{eqnarray}
and then
\begin{eqnarray}\label{M.asy.reg.2}
{\bf P}\{M>x\} &\sim&
\frac{\overline G_2(x)}{\Gamma(\alpha)\Gamma(2-\alpha)}
\quad\mbox{ as } x\to\infty.
\end{eqnarray}
\end{Theorem}

The proofs of Theorems \ref{bounds.gen} and \ref{asy.reg.-}
are given in Section \ref{base.theorems}.
Theorem \ref{asy.reg.-} answers some questions on the behaviour of the
maximums of sums of independent random variables
raised by E.~B.~Dynkin in [\ref{Dynkin}, \S~7].
Some related results for L\'evy processes can be found
in [\ref{KKM}].

Both the tails (\ref{tail.G.1}) and (\ref{tail.G.2}) are lighter than
the integrated tail $\overline{F^I}$ (if the latter exists).

If both the tail $\overline F(t)$ and the function $m(t)$ are regularly
varying at infinity, we can specify the assertion of Theorem
\ref{asy.reg.-} in the following way (the corresponding
calculations are carried out in Section \ref{base.theorems}):

\begin{Corollary}\label{asy.reg.both}
Suppose ${\bf E}\xi^-=\infty$ and the condition {\rm(\ref{finite.M})} holds.
Let $\overline F(t)=t^{-\beta}L^*(t)$ and $m(t)=t^{1-\alpha}L_*(t)$,
where $L^*(t)$ and $L_*(t)$ are functions that are
slowly varying at infinity, $0<\alpha\le1$, $\alpha\le\beta$.
If $\alpha<\beta$, then
\begin{eqnarray}\label{alpha.less.beta}
{\bf P}\{M>x\} &\sim&
\frac{\Gamma(\beta-\alpha)}{\Gamma(\beta)\Gamma(2-\alpha)}
\frac{x\overline F(x)}{m(x)}.
\end{eqnarray}
If  $\alpha=\beta$, then
\begin{eqnarray}\label{alpha.=.beta}
{\bf P}\{M>x\} &\sim&
\frac{1}{\Gamma(\alpha)\Gamma(2-\alpha)}
\int_x^\infty\frac{\overline F(t)}{m(t)}\,dt
\equiv \frac{1}{\Gamma(\alpha)\Gamma(2-\alpha)}
\int_x^\infty\frac{L^*(t)}{tL_*(t)}\,dt.
\end{eqnarray}
\end{Corollary}

\remark
Let $\alpha\in[0,1)$ and $L(x)$ be a slowly varying at infinity
function. Then $m(x) \sim x^{1-\alpha}L(x)$ as $x\to\infty$ if and only if
$F(-x) \sim (1-\alpha)x^{-\alpha}L(x)$
(see [\ref{F}, Chapter XIII, Section 5]).

An asymptotic equivalence like (\ref{alpha.less.beta}) for $\alpha\in(0,1)$,
$\alpha<\beta$ is established in [\ref{Borovkov}, Theorem 4.1]
by other methods and under some additional technical assumptions.
With regard to (\ref{alpha.=.beta}), note that, for any fixed $A>0$,
\begin{eqnarray*}
\int_x^\infty\frac{L^*(t)}{tL_*(t)}\,dt
&\sim& \int_{Ax}^\infty\frac{L^*(t)}{tL_*(t)}\,dt
\quad\mbox{ as } x\to\infty,
\end{eqnarray*}
since, by the Uniform Convergence Theorem for regularly varying
functions (see Theorem 1.5.2 in [\ref{BGT}])
and by Karamata's Theorem (see Proposition 1.5.9b in [\ref{BGT}])
\begin{eqnarray}\label{slow.slow}
\int_x^{Ax}\frac{L^*(t)}{tL_*(t)}\,dt
&\sim& \frac{L^*(x)}{L_*(x)}\,\ln A
=o\Bigl(\int_x^\infty\frac{L^*(t)}{tL_*(t)}\,dt\Bigr).
\end{eqnarray}

Sufficient conditions for the subexponentiality of the distributions
(\ref{tail.G.1}) and (\ref{tail.G.2}) are given in~Section \ref{sufficient.for.G}.
In particular, $G_1$ and $G_2$ are subexponential distributions
if $F$ is either a Pareto, Log-normal or Weibull distribution.
However, in general, the subexponentiality of $F$ only does not imply
that of $G_1$ and $G_2$ (see Section \ref{examples.int}).

The paper is organized as follows. In Sections
\ref{descending.height} and \ref{ascending.height}, we prove some
auxiliary results concerning the first descending and ascending
ladder heights of a random walk.
In Section \ref{base.theorems}, we give the proofs of the
theorems concerning the asymptotics for ${\bf P}\{M>x\}$. Sufficient
conditions for the subexponentiality of (\ref{tail.G.1}) and (\ref{tail.G.2})
may be found in Section \ref{sufficient.for.G}.
Finally, Section \ref{examples.int} is devoted to examples.

\mysection{Asymptotics and bounds for the first descending ladder height\\
in the infinite mean case\label{descending.height}}

Let $\eta_* = \min\{n\ge 1: S_n\le 0\}$ be the first descending
ladder epoch (we put $\min\varnothing=\infty$) and
$\chi_* =-S_{\eta_*}$ be the corresponding descending ladder height.
Since $M$ is finite, $\eta_*$ and $\chi_*$ are proper random variables.
Moreover (see, e.g., Theorem 2.3(c) in [\ref{Aapq}, Chapter VII]),
${\bf E}\eta_*<\infty$ and
\begin{eqnarray}\label{p.eta}
p \equiv {\bf P}\{M=0\} &=& 1/{\bf E}\eta_*.
\end{eqnarray}

For the stopping time $\eta_*$, we have Wald's identity
${\bf E}\chi_*=-{\bf E}\eta_*{\bf E}\xi$, provided
the mean value of $\xi$ is finite and negative
(see Theorem 2(ii) in [\ref{F}, Chapter XII, Section 2]).
In our analysis of the infinite-mean case, the key role
will be played by the following analogue of this identity:

\begin{Lemma}\label{asympt.chi}
Suppose ${\bf E}\xi^-=\infty$ and the condition {\rm(\ref{finite.M})} holds.
Then
\begin{eqnarray}\label{chi.conv}
\frac{{\bf E}\min\,\{\chi_*,x\}}{m(x)} &\to& {\bf E}\eta_*
\quad\mbox{ as } x\to\infty.
\end{eqnarray}
In addition, for any $x\ge0$,
\begin{eqnarray}\label{chi.bound}
{\bf E}\min\,\{\chi_*,x\} &\le& m(x){\bf E}\eta_*.
\end{eqnarray}
\end{Lemma}

\proof. Define the taboo renewal measure on ${\bf R}$
\begin{eqnarray*}
H^*(B) &=& {\bf I}\{0\in B\}+\sum_{n=1}^\infty
{\bf P}\{S_1>0,\ ...,\ S_n>0,\ S_n\in B\}.
\end{eqnarray*}
This measure is finite since $H^*((-\infty,0))=0$ and
\begin{eqnarray}\label{total.H}
H^*([0,\infty)) &=& 1+\sum_{n=1}^\infty
{\bf P}\{S_1>0,\ \ldots,\ S_n>0\}\nonumber\\
&=& 1+\sum_{n=1}^\infty {\bf P}\{\eta_*>n\}
={\bf E}\eta_*<\infty.
\end{eqnarray}
By the total probability formula, for $u\le 0$,
\begin{eqnarray*}
{\bf P}\{-\chi_*\le u\} &=& \int_0^\infty F(u-t)H^*(dt).
\end{eqnarray*}
Therefore,
\begin{eqnarray}\label{repr.for.chi}
\frac{{\bf E}\min\{\chi_*,\,x\}}{m(x)}
&=& \frac{1}{m(x)}\int_0^x {\bf P}\{\chi_*\ge u\}du\nonumber\\
&& \hspace{10mm}= \frac{1}{m(x)}
\int_0^x \int_0^\infty F(-u-t)\,H^*(dt)\,du\nonumber\\
&& \hspace{30mm}= \int_0^\infty \frac{m(x+t)-m(t)}{m(x)} \,H^*(dt).
\end{eqnarray}

For any fixed $z\ge 0$, the function $\min\{z,x\}$ is concave in $x>0$.
Hence, the function $m(x)={\bf E}\min\{\xi^-,x\}$ is concave as well.
In particular, the function $m(x)$ is long-tailed.
Taking into account also that $m(x)\to\infty$ as $x\to\infty$
(since ${\bf E}\xi^-=\infty$),
we deduce the convergence, for any fixed $t\ge 0$,
\begin{eqnarray*}
\frac{m(x+t)-m(t)}{m(x)} &\to& 1 \quad\mbox{ as }x\to\infty.
\end{eqnarray*}
By $m(0)=0$ and by the concavity of $m(x)$,
\begin{eqnarray}\label{incr.m.above}
\frac{m(x+t)-m(t)}{m(x)}
= \frac{m(x+t)-m(t)}{m(x)-m(0)} &\le& 1.
\end{eqnarray}
Applying now the dominated convergence theorem to the finite measure $H^*$,
we obtain the following convergence of the integrals, as $x\to\infty$:
\begin{eqnarray*}
\int_0^\infty \frac{m(x+t)-m(t)}{m(x)} \,H^*(dt)
&\to& \int_0^\infty H^*(dt)=H^*([0,\infty))={\bf E}\eta_*,
\end{eqnarray*}
by (\ref{total.H}). Together with (\ref{repr.for.chi}), this implies
the convergence (\ref{chi.conv}). The inequality (\ref{chi.bound})
follows from (\ref{incr.m.above}) and (\ref{repr.for.chi}).
The proof is complete.

Let $\chi_{*1}$, $\chi_{*2}$, \ldots be independent copies of $\chi_*$.
Define a renewal measure on ${\bf R}^+$
\begin{eqnarray*}
H_*(B) &\equiv& {\bf I}\{0\in B\}
+\sum_{n=1}^\infty {\bf P}\{\chi_{*1}+\cdots+\chi_{*n}\in B\}.
\end{eqnarray*}
If ${\bf E}\xi$ is finite and negative, then
$H_*([0,x])\sim x{\bf E}\chi_*$ as $x\to\infty$,
by the Key Renewal Theorem.
When ${\bf E}\xi$ is infinite, we know only lower and upper estimates
in general:

\begin{Lemma}[{see [\ref{Erickson}, Lemma 1]
or [\ref{BGT}, Section 8.6.3]}]\label{bound.for.H.general}
Without any assumptions, for every $x\ge 0$,
$$
\frac{x}{{\bf E}\min\{\chi_*,x\}} \le H_*([0,x])
\le \frac{2x}{{\bf E}\min\{\chi_*,x\}}.
$$
\end{Lemma}

However, in the regularly varying case, the asymptotic behaviour
of $H_*([0,x])$ is known:

\begin{Lemma}[{see [\ref{Erickson1970}, Theorem 5]}]\label{asym.for.H.general}
If the function ${\bf E}\min\{\chi_*,x\}$ is regularly varying at infinity
with index $1-\alpha$, $\alpha\in[0,1]$, then $H_*([0,x])$ is
regularly varying at infinity with index $\alpha$ and
$$
H_*([0,x]) \sim \frac{1}{\Gamma(1+\alpha)\Gamma(2-\alpha)}\cdot
\frac{x}{{\bf E}\min\{\chi_*,x\}} \quad\mbox{ as }\ x\to\infty.
$$
\end{Lemma}

Using Lemma \ref{asympt.chi} and the equality (\ref{p.eta}),
we obtain from Lemmas \ref{bound.for.H.general}
and \ref{asym.for.H.general} the following corollaries.

\begin{Corollary}\label{cor.bound}
Suppose ${\bf E}\xi^-=\infty$ and the condition {\rm(\ref{finite.M})} holds.
Then
\begin{eqnarray*}
p \le \liminf_{x\to\infty}\frac{H_*([0,x])m(x)}{x}
&\le& \limsup_{x\to\infty}\frac{H_*([0,x])m(x)}{x} \le 2p.
\end{eqnarray*}
\end{Corollary}

\begin{Corollary}\label{cor.regular}
Suppose ${\bf E}\xi^-=\infty$ and the condition {\rm(\ref{finite.M})} holds.
If $m(x)$ is regularly varying at infinity with index $1-\alpha$,
$\alpha\in[0,1]$, then, as $x\to\infty$,
\begin{eqnarray*}
H_*([0,x]) &\sim& \frac{p}{\Gamma(1+\alpha)\Gamma(2-\alpha)}\cdot
\frac{x}{m(x)}.
\end{eqnarray*}
\end{Corollary}

\mysection{Asymptotics and bounds for the first ascending ladder height\\
in the infinite mean case\label{ascending.height}}

Let $\eta^* = \min\{n\ge 1: S_n>0\}$ be the
first ascending ladder epoch
and $\chi^* = S_{\eta^*}$ the corresponding first
ascending ladder height.
Since $M$ is finite a.s., $\eta^*$ and $\chi^*$
are defective random variables, i.e.
${\bf P}\{\eta^*<\infty\}=1-p$ by (\ref{p.eta}).

The starting point in our analysis of the distribution of $\chi^*$
is the following representation (see [\ref{F}, Chapter XII, Section 3]):
\begin{eqnarray}\label{repr.for.chi.+}
{\bf P}\{\chi^*>x\}
&=& \int_0^\infty \overline F(x+t)\,H_*(dt).
\end{eqnarray}

\begin{Lemma}\label{integration.parts}
Suppose ${\bf E}\xi^-=\infty$ and the condition {\rm(\ref{finite.M})} holds.
If the distribution $F$ is long-tailed,
then, for any fixed $T\ge 0$,
\begin{eqnarray*}
{\bf P}\{\chi^*>x\}
&\sim& \int_{x+T}^\infty H_*([0,t-x])\,F(dt)
\quad\mbox{ as }x\to\infty.
\end{eqnarray*}
\end{Lemma}

\proof. Since $F$ is long-tailed and $H_*([0,\infty))=\infty$,
\begin{eqnarray}\label{F.o.int.H}
\overline F(x)
&=& o\Biggl(\int_0^\infty\overline F(x+t)\,H_*(dt)\Biggr)
\quad\mbox{ as }x\to\infty.
\end{eqnarray}
Integration of (\ref{repr.for.chi.+}) by parts gives
\begin{eqnarray}\label{repr.for.chi.+.2}
{\bf P}\{\chi^*>x\}
&=& \overline F(x+t)H_*([0,t])\Big|_0^\infty
+\int_x^\infty H_*([0,t-x])\,d_t F(t).
\end{eqnarray}
Using the upper bound of Corollary \ref{cor.bound}, we obtain,
for sufficiently large $t$,
\begin{eqnarray*}
\overline F(x+t)H_*([0,t])
&\le& \overline F(t)H_*([0,t])
\le 3p\overline F(t)\frac{t}{m(t)}
= 3p\int_t^\infty \frac{t}{m(t)}\,F(ds).
\end{eqnarray*}
Since the function $\frac{x}{m(x)}$ is increasing
(see (\ref{derivative.of.t.mt})),
\begin{eqnarray*}
\overline F(x+t)H_*([0,t])
&\le& 3p\int_t^\infty \frac{s}{m(s)}\,F(ds)\to0
\quad\mbox{ as }t\to\infty,
\end{eqnarray*}
due to condition (\ref{finite.M}).
Substituting this into (\ref{repr.for.chi.+.2}), we arrive
at the equality (recall that $H_*(\{0\})=1$)
\begin{eqnarray*}
{\bf P}\{\chi^*>x\}
&=& -\overline F(x)+\int_x^\infty H_*([0,t-x])\,F(dt).
\end{eqnarray*}
Applying now the relation (\ref{F.o.int.H}),
we deduce the equivalence of the lemma.

In the same way we obtain the following

\begin{Lemma}\label{integration.parts.mt}
Suppose ${\bf E}\xi^-=\infty$ and the condition {\rm(\ref{finite.M})} holds.
If the distribution $F$ is long-tailed,
then, for any fixed $T\ge 0$,
\begin{eqnarray*}
\int_0^\infty \overline F(x+t)\,d\frac{t}{m(t)}
&\sim& \int_{x+T}^\infty \frac{t-x}{m(t-x)}\,F(dt)
\quad\mbox{ as }x\to\infty.
\end{eqnarray*}
\end{Lemma}

\begin{Lemma}\label{chi.+.bound}
Suppose ${\bf E}\xi^-=\infty$ and the condition {\rm(\ref{finite.M})} holds.
If the distribution $F$ is long-tailed, then
\begin{eqnarray*}
p \le \liminf_{x\to\infty}\frac{{\bf P}\{\chi^*>x\}}
{\int_0^\infty \overline F(x+t)\,d\frac{t}{m(t)}}
&\le& \limsup_{x\to\infty}\frac{{\bf P}\{\chi^*>x\}}
{\int_0^\infty \overline F(x+t)\,d\frac{t}{m(t)}} \le 2p.
\end{eqnarray*}
\end{Lemma}

\proof.
Fix $\varepsilon>0$. It follows from Corollary \ref{cor.bound}
that there exists $T>0$ such that, for $t>T$,
\begin{eqnarray*}
(p-\varepsilon)\frac{t}{m(t)} \le
H_*([0,t]) &\le& (2p+\varepsilon)\frac{t}{m(t)}.
\end{eqnarray*}
Applying Lemma \ref{integration.parts}, we obtain,
for $x$ sufficiently large,
\begin{eqnarray*}
(p-2\varepsilon)\int_{x+T}^\infty \frac{t-x}{m(t-x)}\,F(dt)
\le {\bf P}\{\chi^*>x\}
&\le& (2p+2\varepsilon)\int_{x+T}^\infty \frac{t-x}{m(t-x)}\,F(dt).
\end{eqnarray*}
The asymptotic equivalence in Lemma \ref{integration.parts.mt} completes the proof, 
since $\varepsilon>0$ was choosen arbitrary.

Using Corollary \ref{cor.regular} instead of Corollary
\ref{cor.bound}, we may deduce the following

\begin{Lemma}\label{chi.+.regular}
Suppose ${\bf E}\xi^-=\infty$ and the condition {\rm(\ref{finite.M})} holds.
Let the function $m(x)$ be regularly varying at infinity with index
$1-\alpha$, $\alpha\in[0,1]$.
If the distribution $F$ is long-tailed, then
\begin{eqnarray*}
{\bf P}\{\chi^*>x\}
&\sim& \frac{p}{\Gamma(1+\alpha)\Gamma(2-\alpha)}
\int_0^\infty \overline F(x+t)\,d\frac{t}{m(t)}
\quad\mbox{ as }x\to\infty.
\end{eqnarray*}
\end{Lemma}

\proof. It follows from Corollary \ref{cor.regular}
that there exists $T>0$ such that, for $t>T$,
\begin{eqnarray*}
\frac{p-\varepsilon}{\Gamma(1+\alpha)\Gamma(2-\alpha)}
\cdot \frac{t}{m(t)} \le
H_*([0,t]) &\le&
\frac{p+\varepsilon}{\Gamma(1+\alpha)\Gamma(2-\alpha)}
\cdot \frac{t}{m(t)}.
\end{eqnarray*}
By Lemma \ref{integration.parts}, for $x$ sufficiently large,
\begin{eqnarray*}
\frac{p-2\varepsilon}{\Gamma(1+\alpha)\Gamma(2-\alpha)}
\int_{x+T}^\infty \frac{t-x}{m(t-x)}\,F(dt)
&\le& {\bf P}\{\chi^*>x\}\\
&\le& \frac{p+2\varepsilon}{\Gamma(1+\alpha)\Gamma(2-\alpha)}
\int_{x+T}^\infty \frac{t-x}{m(t-x)}\,F(dt).
\end{eqnarray*}
Applying Lemma \ref{integration.parts.mt} completes the proof.

\mysection{The asymptotics and bounds for the distribution tail
of the supremum\label{base.theorems}}

We start with a general theorem which describes the tail behaviour
of the supremum in terms of the renewal measure $H_*$.

\begin{Theorem}\label{base.renewal}
Suppose ${\bf E}\xi^-=\infty$ and the condition {\rm(\ref{finite.M})} holds.
Let the distribution $F$ be long-tailed and the distribution
$G$ with the tail
\begin{eqnarray*}
\overline G(x) &=&
\min\Biggl(1, \int_0^\infty \overline F(x+t)\,d\frac{t}{m(t)}\Biggr)
\end{eqnarray*}
be subexponential. Then, as $x\to\infty$,
\begin{eqnarray*}
{\bf P}\{M>x\} &\sim&
\frac{1}{p}\int_0^\infty \overline F(x+t)\,H_*(dt).
\end{eqnarray*}
\end{Theorem}

\proof.
Consider the distribution $G_H$ with the tail
\begin{eqnarray*}
\overline{G_H}(x) &=&
\min\Biggl(1, \int_0^\infty \overline F(x+t)\,H_*(dt)\Biggr).
\end{eqnarray*}
This distribution is long-tailed because $F$ is long-tailed.
In addition, by Lemma \ref{chi.+.bound}, the tail of $G_H$ is sandwiched
asymptotically between the subexponential tails $p\overline G$ and
$2p\overline G$. Therefore, by the weak equivalence property
(see Theorem 2.1 in [\ref{Kl}] or Lemma 1 in [\ref{AFK}])
the distribution $G_H$ is subexponential as well.

Let us define non-defective random variable
$\widetilde\chi$ with distribution on $(0,\infty)$
\begin{eqnarray*}
{\bf P}\{\widetilde\chi\in B\} &=&
\frac{{\bf P}\{\chi^*\in B\}}{1-p}.
\end{eqnarray*}
The distribution of $\widetilde\chi$ is subexponential.
Let $\widetilde\chi_1$, $\widetilde\chi_2$, \ldots\
be independent copies of the random variable $\widetilde\chi$.
Notice that there exists $i\ge 1$ such that $S_i$ exceeds
a level $x$ if and only if one of the ladder heights
exceeds this level.
Hence, by the formula of total probability we have the equality:
\begin{eqnarray*}
{\bf P}\{M\in B\}
&=& \sum_{n=1}^\infty (1-p)^n p\,
{\bf P}\{\widetilde\chi_1+\cdots+\widetilde\chi_n\in B\}.
\end{eqnarray*}
Since the random variable $\widetilde\chi$ has a subexponential distribution,
we may apply the stopping time theorem
(see, e.g., Lemma 1.8 [\ref{A}, Chapter IX, Section 1]
or Lemma 1.3.5 [\ref{EKM}, Section 1.3.2]) and write
\begin{eqnarray*}
{\bf P}\{M>x\}
&\sim& {\bf P}\{\widetilde\chi>x\}\sum_{n=1}^\infty (1-p)^n p\, n
=\frac{{\bf P}\{\chi^*>x\}}{p}.
\end{eqnarray*}
The proof is complete.

The latter result looks strange in the sense that while the conditions
are expressed in terms of the reference distribution $F$,
the resulting integral is taken with respect to the
renewal measure which is a rather complicated object.
In general, we are unable to write the asymptotics for the integral
$$
\int_0^\infty \overline F(x+t)\,H_*(dt)
$$
in terms of the distribution $F$ itself,
due to the lack of the information about
the asymptotic behaviour of the renewal function $H_*([0,x])$
as $x\to\infty$ in the case of infinite mean.
We may deduce the lower and upper bounds only: combining
the asymptotics in Theorem \ref{base.renewal} and the bounds
in Lemma \ref{chi.+.bound}, we get the assertion of Theorem \ref{bounds.gen}.

To the best of our knowledge, the case when the function $m(x)$
is regularly varying is the only one where the asymptotic behaviour
of $H_*([0,x])$ is known. In this case, combining Theorem
\ref{base.renewal} and Lemma \ref{chi.+.regular},
we obtain the relation (\ref{M.asy.reg.1}) of Theorem \ref{asy.reg.-}.

For $\alpha\in[0,1]$, $tF(-t)=(1-\alpha+o(1))m(t)$ as $t\to\infty$.
Thus, it follows from (\ref{derivative.of.t.mt}) that
\begin{eqnarray*}
\frac{d}{dt}\frac{t}{m(t)} &=& \frac{\alpha+o(1)}{m(t)}
\quad\mbox{ as } t\to\infty.
\end{eqnarray*}
For $\alpha\in(0,1]$, we can apply this result to deduce
(\ref{M.asy.reg.2}) from (\ref{M.asy.reg.1}).

Finally, we prove Corollary \ref{asy.reg.both}.
Notice that the distribution (\ref{tail.G.1}) is subexponential
in this case by Lemma \ref{suff.D.H} from the next section.
We start with the case $0<\alpha\le1$, $\alpha<\beta$.
Fix $\varepsilon>0$ and $A>0$. We have
\begin{eqnarray}\label{0.epsilon.x}
\int_1^{\varepsilon x} \frac{\overline F(x+t)}{m(t)}\,dt
&\le& \overline F(x) \int_1^{\varepsilon x} \frac{1}{m(t)}\,dt
\sim \frac{\overline F(x)}{\alpha}
\frac{\varepsilon x}{m(\varepsilon x)}
\quad\mbox{ as }x\to\infty
\end{eqnarray}
and
\begin{eqnarray}\label{A.x.infty}
\int_{Ax}^\infty \frac{\overline F(x+t)}{m(t)}\,dt
&\le& \int_{Ax}^\infty \frac{\overline F(t)}{m(t)}\,dt
\sim \frac{1}{\beta-\alpha}\frac{Ax\overline F(Ax)}{m(Ax)}
\quad\mbox{ as }x\to\infty.
\end{eqnarray}
Next,
\begin{eqnarray}\label{epsilon.x.A.x}
\int_{\varepsilon x}^{Ax} \overline F(x+t)\frac{1}{m(t)}\,dt
&=& \frac{\overline F(x)}{m(x)}\int_{\varepsilon x}^{Ax}
\frac{\overline F(x+t)}{\overline F(x)}\frac{m(x)}{m(t)}\,dt
\nonumber\\
&=& \frac{x\overline F(x)}{m(x)}\int_\varepsilon^A
\frac{\overline F(x(1+s))}{\overline F(x)}\frac{m(x)}{m(xs)}\,ds
\nonumber\\
&\sim& \frac{x\overline F(x)}{m(x)}\int_\varepsilon^A
\frac{(1+s)^{-\beta}}{s^{1-\alpha}}\,ds
\quad\mbox{ as }x\to\infty,
\end{eqnarray}
since, by the Uniform Convergence Theorem for regularly varying
functions (see Theorem 1.5.2 in [\ref{BGT}])
\begin{eqnarray*}
\frac{\overline F(x(1+s))}{\overline F(x)}\frac{m(x)}{m(xs)}
&\to& \frac{(1+s)^{-\beta}}{s^{1-\alpha}}
\end{eqnarray*}
as $x\to\infty$ uniformly in $s\in[\varepsilon,A]$.
Letting $\varepsilon\to0$ and $A\to\infty$,
we obtain from (\ref{0.epsilon.x}), (\ref{A.x.infty})
and (\ref{epsilon.x.A.x}) that
\begin{eqnarray*}
{\bf P}\{M>x\} &\sim&
\frac{1}{\Gamma(\alpha)\Gamma(2-\alpha)}
\frac{x\overline F(x)}{m(x)}
\int_0^\infty (1+s)^{-\beta}s^{\alpha-1}\,ds
= \frac{B(\beta-\alpha,\alpha)}{\Gamma(\alpha)\Gamma(2-\alpha)}
\frac{x\overline F(x)}{m(x)},
\end{eqnarray*}
which implies (\ref{alpha.less.beta}); here $B$ is the {\it Beta function}.

We now consider the case $0<\alpha\le1$, $\alpha=\beta$.
Fix $A>0$. Now we have
\begin{eqnarray}\label{0.A.x.=}
\int_1^{Ax} \frac{\overline F(x+t)}{m(t)}\,dt
&\le& \overline F(x)\int_1^{Ax} \frac{dt}{m(t)}
\sim \frac{\overline F(x)}{\alpha} \frac{Ax}{m(Ax)}\nonumber\\
&\sim& \frac{A^\alpha}{\alpha} \frac{L^*(x)}{L_*(x)}
=o\Bigl(\int_x^\infty\frac{L^*(t)}{tL_*(t)}\,dt\Bigr),
\end{eqnarray}
since $\alpha=\beta$ and using (\ref{slow.slow}).
Further, for any small $\delta>0$ there exists $A$ sufficiently large
such that $\overline F(x+t)\ge(1-\delta)\overline F(t)$ for any
$t\ge Ax$. Then
\begin{eqnarray}\label{A.x.infty.=.1}
\int_{Ax}^\infty \frac{\overline F(x+t)}{m(t)}\,dt
&\ge& (1-\delta)\int_{Ax}^\infty \frac{\overline F(t)}{m(t)}\,dt
= (1-\delta)\int_{Ax}^\infty\frac{L^*(t)}{tL_*(t)}\,dt.
\end{eqnarray}
On the other hand,
\begin{eqnarray}\label{A.x.infty.=.2}
\int_{Ax}^\infty \frac{\overline F(x+t)}{m(t)}\,dt
&\le& \int_{Ax}^\infty \frac{\overline F(t)}{m(t)}\,dt
= \int_{Ax}^\infty\frac{L^*(t)}{tL_*(t)}\,dt.
\end{eqnarray}
The relations (\ref{0.A.x.=}), (\ref{A.x.infty.=.1}),
(\ref{A.x.infty.=.2}) and (\ref{slow.slow}) imply
(\ref{alpha.=.beta}). Corollary \ref{asy.reg.both} is proved.

\mysection{Sufficient conditions for the integrated weighted tail\\
distribution to be subexponential\label{sufficient.for.G}}

In this Section, we present sufficient conditions for the
subexponentiality of the distributions (\ref{tail.G.1}) and (\ref{tail.G.2}).
We consider an even more general problem:
Let $F$ be a distribution on ${\bf R}^+$ and $H$ a non-negative
measure on ${\bf R}^+$ such that
\begin{equation}\label{F.H.finite}
\int_0^\infty \overline F(t)\,H(dt)\quad\mbox{is finite}.
\end{equation}
In this case we can define the distribution $G_H$
on ${\bf R}^+$ with the tail
\begin{eqnarray}\label{integrated.weighted.tail}
\overline G_H(x) &\equiv& \min\Biggl(1,\int_0^\infty
\overline F(x+t)\,H(dt)\Biggr),\quad x\ge0.
\end{eqnarray}
We can formulate the following question:
what type of conditions on $F$ imply the subexponentiality of $G_H$?

First, recall that if $F$ is long-tailed, then $G_H$ is long-tailed as well.

\definition{}
The distribution $F$ on ${\bf R}^+$ is called dominated varying
($F\in\mbox{\scal D}$) iff, for some $c>0$,
$\overline F(2x)\ge c\overline F(x)$ for any $x$.

It is known that $(\mbox{\scal L}\cap\mbox{\scal D})\subset\mbox{\scal S}$.  Also, it is known that if $F\in\mbox{\scal D}$,
then $F^I\in\mbox{\scal L}\cap\mbox{\scal D}$,
but the converse is not true in general (see [\ref{Kl}, Section 4]).

\begin{Lemma}\label{suff.D.H}
If $F\in\mbox{\scal D}\,\cap\mbox{\scal L}$,
then $G_H\in\mbox{\scal D}\,\cap\mbox{\scal L}$ and, therefore,
$G_H\in\mbox{\scal S}$.
\end{Lemma}

\proof. This result follows from the inequalities:
\begin{eqnarray*}
\int_0^\infty \overline F(2x+t)\,H(dt)
&\ge& c\int_0^\infty \overline F(x+t/2)\,H(dt)
\ge c\int_0^\infty \overline F(x+t)\,H(dt).
\end{eqnarray*}

\definition{}
The distribution $F$ on ${\bf R}^+$ with finite mean $m$ belongs
to the class $\mbox{\scal S}\,^*$ if
\begin{eqnarray*}
\int_0^x \overline F(x-y)\overline F(y)dy
&\sim& 2m \overline F(x)\quad\mbox{ as }x\to\infty.
\end{eqnarray*}
It is known (see [\ref{Kl}]) that
\begin{equation} \label{conj1}
F\in\mbox{\scal S}\,^* \quad \mbox{implies} \quad
F\in\mbox{\scal S}\,
\quad \mbox{and} \quad  F^I\in\mbox{\scal S}.
\end{equation}
It turns out that the following more general conclusion holds.
For any $b>0$, define the class $\mbox{\scal H}_b$ of all
non-negative measures $H$ on ${\bf R}^+$ such that
\begin{eqnarray*}
\sup_t H((t,t+1]) &\le& b.
\end{eqnarray*}

\begin{Lemma} \label{L9}
Let $F\in\mbox{\scal S}\,^*$ and $H\in\mbox{\scal H}_b$, $b\in(0,\infty)$.
Then $G_H\in\mbox{\scal S}$. Moreover,
\begin{eqnarray*}
\overline{G_H*G_H}(x) &\sim& 2\overline G_H(x)
\end{eqnarray*}
as $x\to\infty$ uniformly in $H\in\mbox{\scal H}_b$.
\end{Lemma}

\remark
Here are four examples of such measures $H$:
(i) if $H(B)={\bf I}\{0\in B\}$, then $G_H=F$;
(ii) if $H(dt)=dt$ is Lebesgue measure on ${\bf R}^+$, then $G_H=F^I$;
(iii) if $H$ is the renewal measure $H_*$, then $G_H$ is the
distribution of the first ascending ladder height $\chi^*$;
(iv) if $H([0,x])=x/m(x)$, then $G_H$ is $G_1$ from (\ref{tail.G.1}).

\remark
It is natural to consider the following two questions:

(i) may the assumption $F\in\mbox{\scal S}\,^*$
of Lemma \ref{L9} be weakened to $F\in\mbox{\scal S}\,$?
In the case of Lebesgue measure $H$, i.e. when $G_H=F^I$,
this question is raised in [\ref{EKM}, Section 1.4.2].

(ii) is the converse of (\ref{conj1}) also true?

In the next Section, we show (by examples) that
the answers to both these questions are negative.

\proof\ \ of Lemma \ref{L9}. Since $G_H$ is long-tailed
uniformly in $H\in\mbox{\scal H}_b$, it is sufficient to show that
\begin{eqnarray}\label{G.suff.S.*}
\lim_{A\to\infty}\limsup_{x\to\infty} \sup_{H\in\mbox{\scall H}_b}
\frac{1}{\overline G_H(x)} \int_A^{x-A} \overline G_H(x-y) G_H(dy)
&=& 0,
\end{eqnarray}
see, e.g., Proposition 2 in [\ref{AFK}].

The mean value of $F$ is finite. Thus,
$\overline F(t)H((0,t])=o(1/t)O(t)\to0$ as $t\to\infty$ and
integration by parts yields, for $x$ large enough,
$$
\overline G_H(x) = \int_x^\infty H((0,t-x])F(dt).
$$
Hence,
\begin{eqnarray*}
G_H((x,x+1]) &=& \int_x^\infty H((t-x-1,t-x])\,F(dt)
\le b\overline F(x).
\end{eqnarray*}
In addition, $G_H$ is long-tailed.
Therefore, (\ref{G.suff.S.*}) holds if and only if
\begin{eqnarray}\label{G.suff.S.*.2}
\lim_{A\to\infty}\limsup_{x\to\infty} \sup_{H\in\mbox{\scall H}_b}
\frac{1}{\overline G_H(x)}
\int_A^{x-A} \overline G_H(x-y) \overline F(y)\,dy
&=& 0.
\end{eqnarray}

Fix $\varepsilon>0$. Since $F\in\mbox{\scal S}\,^*$, there exist
$x_0$ and $A$ such that, for all $x\ge x_0$,
$$
\int_A^{x-A} \overline F(x-u)\overline F(u)\,du
\le \varepsilon\overline F(x).
$$
Then, for $x\ge x_0$,
\begin{eqnarray*}
\int_A^{x-A} \overline G_H(x-y) \overline F(y)\,dy
&=& \int_A^{x-A} \Biggl(\int_0^\infty\overline F(x+t-y)\,H(dt)\Biggr)
\overline F(y)dy\\
&\le& \int_0^\infty \Biggl(\int_A^{x+t-A} \overline F(x+t-y) \overline F(y)
dy\Biggr) H(dt)\\
&\le& \varepsilon \int_0^\infty \overline F(x+t)H(dt)
=\varepsilon\overline G_H(x).
\end{eqnarray*}
Letting $\varepsilon$ to $0$, we get (\ref{G.suff.S.*.2}).
The proof is complete.

Lemma \ref{L9} also implies that any $\mbox{\scal S}\,^*$-distribution 
$F$ is strongly subexponential in the sense of [\ref{K2}], i.e., 
\begin{eqnarray*}
\overline{F_h * F_h}(x) &\sim& 2\overline F_h(x)
\end{eqnarray*}
as $x\to\infty$ uniformly in $h\in[1,\infty)$, 
where the distribution $F_h$ is defined as follows:
\begin{eqnarray*}
\overline{F_h}(x) &\equiv& 
\min\Bigl(1,\int_x^{x+h}\overline F(t)dt\Bigr),\quad x>0.
\end{eqnarray*}
Due to the main result of [\ref{K2}], 
it allows us to formulate the following

\begin{Corollary}
Let $\xi$ have finite negative mean value and 
its distribution $F$ be from $\mbox{\scal S}\,^*$. Then
\begin{eqnarray*}
{\bf P}\Bigl\{\max_{0\le k\le n}S_k>x\Bigr\}
&\sim& \frac{1}{|{\bf E}\xi|}
\int_x^{x+n|{\bf E}\xi|} \overline F(t)dt
\end{eqnarray*}
as $x\to\infty$ uniformly in $n\ge 1$.
\end{Corollary}

\mysection{Examples\label{examples.int}}

In this Section, we give an example of $F\in\mbox{\scal S}$ with
finite mean such that $F^I\notin\mbox{\scal S}$.
In fact, we provide a more general example:
for any fixed $\alpha\in[0,1)$, we construct a subexponential
distribution $F$ with finite mean such that
the distribution $G_\alpha$ with the tail
\begin{equation}\label{Ga}
\overline{G_\alpha}(x)
= \min\Bigl(1,\int_1^\infty\frac{\overline F(x+y)}{y^\alpha}
\,dy
\Bigr)
\end{equation}
is not subexponential. In particular, when $\alpha=0$,
$F^I $ does not belong to $\mbox{\scal S}$.

In our second example, we show that two conditions $F\in\mbox{\scal S}$
and $F^I\in\mbox{\scal S}$ taken together do not imply that
$F\in\mbox{\scal S}\,^*$.

Both examples are based on the following construction.

Define two increasing sequences of positive numbers,
namely $\{t_n\}$ and $\{R_n\}$, such that, as $n\to\infty$,
\begin{eqnarray}\label{tn}
t_n = o(t_{n+1}),\\
\label{Rn}
R_{n+1}-R_n \to \infty.
\end{eqnarray}
Define the {\it hazard function} $R(x)\equiv -\ln\overline F(x)$ as
$$
R(x)=R_n+r_n(x-t_n)\ \mbox{ for }\ t_n\le x\le t_{n+1},
$$
where
\begin{equation}\label{rn}
r_n=\frac{R_{n+1}-R_n}{t_{n+1}-t_n}
\sim \frac{R_{n+1}-R_n}{t_{n+1}}.
\end{equation}
by (\ref{tn}) and (\ref{Rn}).
In other words, the {\it hazard rate} $r(x)\equiv R'(x)$
is defined as $r(x)=r_n$ for $x\in(t_n,t_{n+1}]$,
where $r_n$ is given by (\ref{rn}).

Note that
$$
J_n \equiv \int_{t_n}^{t_{n+1}} \overline F(u) du
= \int_{t_n}^{t_{n+1}} e^{-R(u)} du
= \frac{e^{-R_n}-e^{-R_{n+1}}}{r_n}
< \frac{e^{-R_n}}{r_n},
$$
and that the mean value of $F$ is finite provided
\begin{equation}\label{mean}
\sum_n \frac{e^{-R_n}}{r_n} < \infty .
\end{equation}
We assume that (\ref{mean}) holds. Assume also that
\begin{equation} \label{r0}
r_{n+1}=o(r_n)\ \mbox{ and }\ r_nt_n\to 0\ \mbox{ as }n\to\infty.
\end{equation}
It follows from (\ref{r0}) that $r_kt_n\to 0$ as $n\to\infty$
uniformly in $k\ge n$. However,
\begin{eqnarray}\label{r.n.t.n1}
r_nt_{n+1} &\sim& R_{n+1}-R_n \to \infty
\end{eqnarray}
from (\ref{rn}) and (\ref{Rn}). It follows from (\ref{r0}) that
$r(x)$ decreases eventually to 0, and we can
apply the following results:

\begin{Proposition}[see Corollary 3.8 and Theorem 3.6 in {[\ref{Kl}]}]
If the hazard rate exists and is eventually decreasing
to $0$, then $F\in\mbox{\scal L}$ and

{\rm(i)} $F\in\mbox{\scal S}$ if and only if
\begin{equation}\label{eqv1}
\lim_{x\to\infty} r(x) \int_0^x e^{yr(x)} \overline F(y) dy = 0.
\end{equation}

{\rm(ii)} $F\in\mbox{\scal S}\,^*$ if and only if $F$ has finite mean and
\begin{equation} \label{eqv2}
\lim_{x\to\infty}\int_0^x e^{yr(x)} \overline F(y)\,dy =
\int_0^\infty \overline F(y)\,dy.
\end{equation}
\end{Proposition}

Note that (\ref{eqv2}) is equivalent to
\begin{equation} \label{eqv3}
\lim_{t\to\infty}\lim_{x\to\infty}
\int_t^x e^{ yr(x)} \overline F(y)\,dy = 0.
\end{equation}
Put
$$
I_{n,k}=\int_{t_k}^{t_{k+1}} e^{yr_n}\overline F(y)\,dy
\quad\mbox{ and }\quad
I_n = \int_1^{t_{n+1}} e^{yr_n} \overline F(y)\,dy
=\sum_{k=1}^n I_{n,k}.
$$
In our case, (\ref{eqv1}) holds if and only if
\begin{equation} \label{rI}
r_n I_n \to 0\ \mbox{ as }n\to\infty.
\end{equation}
The relation (\ref{eqv3}) fails, in particular, if
\begin{equation}\label{fail}
\liminf_{n\to\infty}I_{n,n} > 0.
\end{equation}

From (\ref{r0}),
$$
I_{n,k} = e^{-R_k+r_kt_k}
\frac{e^{(r_n-r_k)t_k}-e^{(r_n-r_k)t_{k+1}}}{r_k-r_n}
\le \frac{e^{-R_k+r_nt_k}}{r_k-r_n}
\sim \frac{e^{-R_k}}{r_k}
$$
as $k\to\infty$ uniformly in $n\ge k+1$. Thus, for some $C<\infty$,
\begin{equation}\label{ink}
\sum_{k=1}^{n-1} I_{n,k} \le
C \sum_{k=1}^{n-1} \frac{e^{-R_k}}{r_k}
\le C \sum_{k=1}^\infty \frac{e^{-R_k}}{r_k} < \infty
\end{equation}
if (\ref{mean}) holds. Further,
\begin{equation}\label{inn}
r_nI_{n,n} = r_n(t_{n+1}-t_n) e^{-R_n+r_nt_n}
\sim (R_{n+1}-R_n) e^{-R_n}
\end{equation}
from (\ref{rn}) and (\ref{r0}). Thus, $r_nI_{n,n} \to 0$ if
\begin{equation}\label{cnR}
R_{n+1}=o(e^{R_n}).
\end{equation}
Hence, under the conditions (\ref{tn}), (\ref{Rn}), (\ref{mean}),
(\ref{r0}), and (\ref{cnR}), $F$ is a subexponential distribution
with finite mean.

We now turn to the examples.

\begin{Example}
Fix $\alpha\in[0,1)$ and put $R_{n+1}=e^{\gamma R_n}$,
where the constant $\gamma=\gamma(\alpha)\in(0,1)$ will be specified later.
If we take $R_1=R_1(\gamma)$ sufficiently large, then the sequence
$R_n$ will be increasing and, moreover, $R_{n+1}/R_n\to\infty$.
Put $t_{n+1}=e^{2\gamma R_n}=R_{n+1}^2$;
the condition (\ref{tn}) is satisfied. We have
$$
\overline F(t_n) = e^{-\sqrt{t_n}}.
$$
We also have $r_n\sim R_{n+1}/t_{n+1}=e^{-\gamma R_n}$.
\end{Example}

The condition (\ref{mean}) is valid since $J_n \sim e^{-R_n(1-\gamma)}$.
The condition (\ref{r0}) holds since
$r_{n+1}/r_n \sim e^{-\gamma(R_{n+1}-R_n)}$ and
$r_nt_n \sim e^{-\gamma R_n + 2\gamma R_{n-1}} = e^{-\gamma R_n+o(R_n)}$.
Finally, the condition (\ref{cnR}) follows since $R_{n+1}e^{-R_n} = e^{-(1-\gamma)R_n}$.
Hence, $F$ has a finite mean and is subexponential.

Take now the distribution $G_\alpha$ defined in (\ref{Ga})
and estimate its density. For $x\in(t_n,t_{n+1}-1]$,
\begin{eqnarray}\label{deri}
G_\alpha'(x) &=& -\int_1^\infty\frac{\overline F'(x+y)}{y^\alpha}\,dy
= \int_1^\infty\frac{r(x+y)\overline F(x+y)}{y^\alpha}\,dy\nonumber\\
&\ge& \int_1^{t_{n+1}-x}\frac{r(x+y)\overline F(x+y)}{y^\alpha}\,dy
= r_n V_n(x),
\end{eqnarray}
where
$$
V_n(x) = \int_1^{t_{n+1}-x} \overline F(x+y) y^{-\alpha} dy.
$$
We also have that, for any $x<t_{n+1}-1$,
\begin{eqnarray}\label{deri.2}
\overline G_\alpha(x) &\ge& V_n(x).
\end{eqnarray}

For any $x\in(t_n,t_{n+1}-1]$,
$$
\overline{G_\alpha}(x) =
\left(\int_{x+1}^{t_{n+1}}+\sum_{k=n+1}^\infty
\int_{t_k}^{t_{k+1}} \right)
\overline F(y) (y-x)^{-\alpha} dy
\equiv V_n(x) + \sum_{k=n+1}^\infty W_k(x).
$$
For $x\in(t_n,t_{n+1}]$ and $k\ge n+1$, by (\ref{r0}) and (\ref{r.n.t.n1}),
\begin{eqnarray*}
W_k(x) &=& e^{-R_k+r_kt_k}
\int_{t_k}^{t_{k+1}} e^{-r_ky} (y-x)^{-\alpha} dy\\
&\sim& e^{-R_k} \int_{t_k}^{t_{k+1}}
e^{-r_ky} (y-x)^{-\alpha} dy
= \frac{e^{-R_k}}{r_k^{1-\alpha}}
\int_{r_kt_k}^{r_kt_{k+1}} e^{-y}(y-r_kx)^{-\alpha} dy\\
&\sim& \frac{e^{-R_k}}{r_k^{1-\alpha}}
\int_{0}^\infty e^{-y}y^{-\alpha} dy
=  \frac{e^{-R_k}}{r_k^{1-\alpha}}\Gamma(1-\alpha)
\quad\mbox{as }k\to\infty.
\end{eqnarray*}
Similarly, for $x\in(t_n,t_{n+1}/2]$,
\begin{eqnarray*}
V_n(x) &=& e^{-R_n+r_nt_n-r_nx}
\int_1^{t_{n+1}-x} e^{-r_ny} y^{-\alpha} dy\\
&\sim& \frac{e^{-R_n-r_nx}}{r_n^{1-\alpha}}
\int_{r_n}^{r_n(t_{n+1}-x)} e^{-z}z^{-\alpha} dz
\sim \frac{e^{-R_n-r_nx}}{r_n^{1-\alpha}} \Gamma(1-\alpha)
\quad\mbox{as }n\to\infty.
\end{eqnarray*}

Since $\gamma<1$,
$$
\frac{W_{k+1}(x)}{W_k(x)} \sim
\left(\frac{r_k}{r_{k+1}}\right)^{1-\alpha} e^{-R_{k+1}+R_k}
\sim e^{[\gamma(1-\alpha )-1](R_{k+1}-R_k)} \to 0.
$$
Take any integer $l\ge 2$ such that $\gamma(1-\alpha)<(l-1)/l$.
Then, as $n\to\infty$,
$$
\frac{W_{n+1}(x)}{V_n(t_{n+1}/l)} \sim
\left(\frac{r_n}{r_{n+1}}\right)^{1-\alpha}
e^{-R_{n+1}+R_n+r_nt_{n+1}/l}
=e^{[\gamma(1-\alpha)-(l-1)/l]R_{n+1}+o(R_{n+1})}\to 0.
$$
Therefore,
$$
\overline{G_\alpha}(t_{n+1}/l) \sim V_n(t_{n+1}/l)
\sim \frac{e^{-R_n-r_nt_{n+1}/l}}{r_n^{1-\alpha}} \Gamma(1-\alpha)
\ \mbox{ as }n\to\infty.
$$
On the other hand, by (\ref{deri}), for $n$ sufficiently large,
\begin{eqnarray*}
\overline{G_\alpha*G_\alpha}(t_{n+1}/l)
&\ge& \int_{t_n}^{t_{n+1}/2l}
\overline G_\alpha(t_{n+1}/l-y)\,G_\alpha(dy)\\
&\ge& r_n\int_{t_n}^{t_{n+1}/2l}
\overline G_\alpha(t_{n+1}/l-y)V_n(y)\,dy.
\end{eqnarray*}
Applying now (\ref{deri.2}), we get
\begin{eqnarray*}
\overline{G_\alpha*G_\alpha}(t_{n+1}/l)
&\ge& r_n\int_{t_n}^{t_{n+1}/2l}
V_n(t_{n+1}/l-y)V_n(y)\,dy\\
&\sim& \frac{\Gamma^2(1-\alpha)}{r_n^{1-2\alpha}}
\int_{t_n}^{t_{n+1}/2l}
e^{-R_n-r_n(t_{n+1}/l-y)-R_n-r_ny}\,dy\\
&\sim& \frac{\Gamma^2(1-\alpha)}{r_n^{1-2\alpha}}
\frac{t_{n+1}}{2l}e^{-2R_n-r_nt_{n+1}/l}.
\end{eqnarray*}
Then the ratio
$$
\frac{\overline{G_\alpha*G_\alpha}(t_{n+1}/l)}
{\overline{G_\alpha}(t_{n+1}/l)}
$$
is asymptotically not less than
$$
\frac{\Gamma (1-\alpha )}{2l}r_n^{\alpha}  t_{n+1}e^{-R_n}
\sim \frac{\Gamma (1-\alpha )}{2l} e^{R_n(-\gamma\alpha+2\gamma-1)}
\to\infty
$$
as $n\to\infty$ provided $\gamma(2-\alpha)>1$. Thus, for any
$\gamma\in(1/(2-\alpha),1)$, $F\in\mbox{\scal S}$ and has finite mean,
but $G_{\alpha}\notin\mbox{\scal S}$.

\begin{Example}
For $\gamma>2$, take $R_n=n^\gamma$
and $t_{n+1}=e^{R_n}=e^{n^\gamma}$. Then
$$
\overline F(t_{n})=t_n^{-\bigl(\frac{n}{n-1}\bigr)^\gamma}.
$$
Conditions (\ref{tn}), (\ref{Rn}), and (\ref{r0}) are satisfied,
$r_n \sim \gamma n^{\gamma-1}/t_{n+1}$,
and (\ref{mean}) holds. Further, (\ref{cnR}) holds too.
Hence, $F\in\mbox{\scal S}$.
\end{Example}

On the other hand, for $x\in(t_n,t_{n+1}]$,
$$
\overline{F^I}(x) \le \sum_{k=n}^\infty J_k
\ \mbox{ and }\ \overline{F^I}(2x) \ge \sum_{k=n+2}^\infty J_k.
$$
In addition, $J_k \sim 1/\gamma k^{\gamma-1}$ as $k\to\infty$.
Thus, $\overline{F^I}(2x) \sim \overline{F^I}(x)$ as $x\to\infty$
and the function $\overline{F^I}(x)$ is slowly varying at infinity.
Hence, $F^I\in\mbox{\scal S}$.

However, from (\ref{inn}) and (\ref{rn}),
$$
I_{n,n}\sim t_{n+1}e^{-R_n}=1
$$
and so, from (\ref{fail}), $F$ cannot belong to $\mbox{\scal S}\,^*$.

\section*{\normalsize References}
\small

\newcounter{bibcoun}
\begin{list}{[\arabic{bibcoun}]}{\usecounter{bibcoun}\itemsep=0pt}
\small

\item\label{Aapq}
S. Asmussen,
{\em Applied Probability and Queues}
(John Wiley \& Sons, Chichester, 1987).

\item\label{A}
S. Asmussen,
{\em Ruin Probabilities}
(World Scientific, Singapore, 2000).

\item\label{AFK}
S.~Asmussen, S.~Foss and D.~Korshunov,
Asymptotics for sums of random variables with
local subexponential behaviour,
J. Theoret. Probab. 16 (2003) 489--518.

\item\label{BGT}
N.~H.~Bingham, C.~M.~Goldie and J.~L.~Teugels,
{\em Regular Variation}
(Cambridge University Press, 1987).

\item\label{Borovkov}
A.~A.~Borovkov,
Large deviations probabilities for random walks
in the absence of finite expectations of jumps,
Probab. Theory Relat. Fields 125 (2003) 421--446.

\item\label{Ch}
V.~P.~Chistyakov,
A theorem on sums of independent random positive
variables and its applications to branching processes,
Theory Probab. Appl. 9 (1964) 710--718.

\item\label{Dynkin}
E.~B.~Dynkin,
Some limit theorems for sums of independent random
variables with infinite mathematical expectations,
Select. Transl. Math. Statist. Probability 1 (1961) 171--189.

\item\label{EKM}
P. Embrechts, C. Kl\"uppelberg and T. Mikosch,
{\em Modelling Extremal Events}
(Springer, Berlin, 1997).

\item\label{EO}
P.~Embrechts and E.~Omey,
A property of longtailed distributions,
J.~Appl. Prob. 21 (1984) 80--87.

\item\label{EV}
P.~Embrechts and N.~Veraverbeke,
Estimates for the probability of ruin with special
emphasis on the possibility of large claims,
Insurance: Mathematics \& Economics 1 (1982) 55--72.

\item\label{Erickson1970}
K.~B.~Erickson,
Strong renewal theorems with infinite mean,
Transactions of the American Mathematical Society
151 (1970) 263--291.

\item\label{Erickson}
K.~B.~Erickson,
The strong law of large numbers when the mean is undefined,
Transactions of the American Mathematical Society
185 (1973) 371--381.

\item\label{F}
W.~Feller,
{\em An Introduction to Probability Theory and Its Applications},
Vol.~2, 2nd ed. (Wiley, New York, 1971).

\item\label{Kl}
C.~Kl\"uppelberg,
Subexponential distributions and integrated tails,
J. Appl. Probab. 25 (1988) 132--141.

\item\label{KKM}
C.~Kl\"uppelberg, A.~E.~Kyprianou and R.~A.~Maller,
Ruin probabilities and overshoots for general L\'evy
insurance risk processes,
submitted for publication (2003).

\item\label{K}
D.~Korshunov,
On distribution tail of the maximum of a random walk,
Stochastic Process. Appl. 72 (1997) 97--103.

\item\label{K2}
D.~A.~Korshunov,
Large-deviation probabilities for maxima of sums of
independent random variables with negative mean and
subexponential distribution,
Theory Probab. Appl. 46 (2002) 355--366.

\item\label{T}
J.~L.~Teugels,
The class of subexponential distributions,
Ann. Probab. 3 (1975) 1000--1011.

\item\label{V}
N.~Veraverbeke,
Asymptotic behavior of Wiener--Hopf factors of a random walk,
Stochastic Process. Appl. 5 (1977) 27--37.

\end{list}

\end{document}